\documentclass{siamart251216}

\usepackage{amssymb}
\usepackage{booktabs}

\IfFileExists{supplement.aux}
    {\externaldocument[supp-]{supplement}[supplement.pdf]}
    {\IfFileExists{.latex-build/supplement-preview/supplement.aux}
        {\externaldocument[supp-]{.latex-build/supplement-preview/supplement}[supplement.pdf]}
        {}}
\newcommand{\suppsection}[2]{%
    \ifcsname r@supp-#1\endcsname
        Section~\ref{supp-#1} of the supplement%
    \else
        the supplement's section \emph{#2}%
    \fi}

\headers{Mesh-uniform residual convergence}{P. Jiang and T. Lu}
\title{Mesh-uniform residual convergence of ADMM for slope-selection epitaxy and the phase-field crystal equation%
\thanks{This work was partially supported by the Fundamental and Interdisciplinary Disciplines Breakthrough Plan of the Ministry of Education of China (Grant No. JYB2025XDXM204).}}
\author{Peng Jiang%
\thanks{School of Mathematical Sciences, Peking University, Beijing, China (\email{pengjiang@stu.pku.edu.cn}).}
\and Tiao Lu%
\thanks{Corresponding author. School of Mathematical Sciences and HEDPS, CAPT, LMAM, Peking University, Beijing, China (\email{tlu@pku.edu.cn}).}}

\begin{document}
\maketitle

\begin{abstract}
We study the number of alternating direction method of multipliers (ADMM) sweeps needed to reach a prescribed tolerance in the projected optimality residual of the original convex-splitting time-step objective. For unforced backward Euler and stabilized second-order backward differentiation formula (BDF2) steps of slope-selection epitaxy and the standard phase-field crystal equation, we obtain a geometric bound for the unnormalized residual in the unpreconditioned discrete \(L^2\) norm, uniform in the mesh size and in \(0<\tau\le\tau_0\). The single-step result assumes bounded histories and initialization as specified below and exact subproblem solves. Fourier smoothing and a pointwise invariant bound control the quartic proximal orbit. A discrete first-difference estimate transfers the Douglas--Rachford contraction to the epitaxy residual. Along exact, unforced trajectories with prescribed startup and restarts the bounds are uniform in the mesh and time level, also in \(\tau\) for the crystal model and at fixed \(\tau\) for epitaxy. The bound allows synchronized relaxation as in ADMM-Search. Smooth-data refinement tests illustrate the sweep counts.

\end{abstract}

\begin{keywords}
Alternating direction method of multipliers, mesh-uniform convergence, projected optimality residual,
convex splitting, slope-selection epitaxy, phase-field crystal equation
\end{keywords}

\begin{AMS}
65M22, 65M12, 90C25
\end{AMS}

\section{Introduction}
\label{sec:introduction}

Convex splitting gives energy-stable discretizations of gradient flows, often with a convex minimization problem at each time step \cite{Eyre1998,ShenEtAl2012}. We consider unforced backward Euler (BE) and stabilized second-order backward differentiation formula (BDF2) steps for slope-selection (SS) molecular beam epitaxy (MBE) and the standard phase-field crystal (PFC) equation. Can the alternating direction method of multipliers (ADMM) reduce the projected gradient of the time-step objective to a prescribed tolerance in a mesh-independent number of sweeps?

For SS, Feng et~al.~\cite{FengSalgadoWangWise2017} prove geometric convergence of preconditioned steepest descent (PSD) for regularized \(p\)-Laplacian problems. Feng et~al.~\cite{FengEtAl2018} develop a stabilized SS-BDF2 scheme and discuss PSD and preconditioned nonlinear conjugate gradient (PNCG) solvers. Their relative residual histories show nearly mesh-independent decay under grid refinement at fixed physical parameters and step size \cite[Section~6.1, Fig.~1]{FengEtAl2018}. Those SS discretizations use vertex gradients, whereas ours lifts face differences to the cells. For standard PFC, existing work establishes solvability, energy stability, and convergence of finite-difference convex splitting \cite{WiseWangLowengrub2009}, and develops nonlinear multigrid \cite{HuWiseWangLowengrub2009} and preconditioned acceleration for Fourier collocation \cite{ParkSalgadoWise2023}. A related accelerated-gradient analysis estimates the residual at the extrapolated iterates in the dual norm induced by their preconditioner \cite[Corollary~5.9]{ParkSalgadoWise2021}.

The Douglas--Rachford (DR) state lies in the split-variable space, and the projected gradient lies in the tangent space of the admissible set. In both models the quartic Hessian vanishes at zero and is unbounded. The standard reflected-proximal contraction estimate based on a positive global lower curvature bound and a finite global upper curvature bound does not apply. Fourier smoothing controls the recovered fields pointwise, and the quartic proximal update preserves a bounded range.

\begin{enumerate}
    \item We formulate a common ADMM solver for a class of convex-split gradient-flow time steps. On a fixed grid, we derive an explicit upper bound on the Q-linear contraction factor of the DR state when the convex potential has a locally Lipschitz gradient (Theorem~\ref{thm:admm-linear}).
    \item Our main result, Theorem~\ref{thm:mesh-uniform-admm}, gives a mesh- and step-size-uniform geometric bound for the unnormalized projected optimality residual of a single SS or standard PFC time step in the unpreconditioned discrete \(L^2\) norm, under the stated history and first-sweep conditions and with exact subproblem solves. For SS, a first-difference estimate before applying \(B_h^*\) avoids an inverse mesh scale.
    \item For exact, unforced trajectories with the prescribed startup and restarts, Corollary~\ref{cor:mesh-uniform-outer-trajectories} extends the residual bound to all time levels, including startup. The bounds are uniform in the mesh and time level for both models. They are also uniform in the time step for PFC and hold at each fixed time step for SS.
\end{enumerate}

We use the classical ADMM--DR equivalence and fixed-point error bounds \cite{EcksteinBertsekas1992,LiuYuanZengZhang2018,PenaVeraZuluaga2021}. The residual inherits the state's geometric factor, yielding a sweep bound for a fixed tolerance. The analysis covers the synchronized relaxation of ADMM-Search. Section~\ref{sec:numerical-results} illustrates the mesh dependence of the sweep counts under the normalized stopping criterion. The supplement gives further experiments and implementation details.

\section{Discrete time-step problems and model operators}
\label{sec:discretization}

\subsection{Gradient flows and grid spaces}
\label{sec:ccfd}
We consider a scalar field \(u\) on the periodic box \(\Omega=\prod_{j=1}^d(0,\ell_j)\). With \(\langle v,w\rangle=\int_\Omega v\cdot w\,\mathrm dx\), the continuous split energy and its gradient flow take the form
\begin{gather}
 E(u)=\tfrac12\langle Au,u\rangle+\int_\Omega F(Bu)\,\mathrm dx
                  -\tfrac12\langle Cu,u\rangle,\label{eq:continuous-gradient-flow-energy}\\
 \partial_t u=-Lw,\qquad w=Au+B^*f(Bu)-Cu .\label{eq:continuous-gradient-flow}
\end{gather}
Here \(A,C,L\) are self-adjoint and positive semidefinite, \(B\) is linear, and \(F:\mathbb R^m\to\mathbb R\) is convex and differentiable, with \(f=\nabla F\). The convex energy consists of the positive quadratic term and the nonlinear density. The concave contribution is the negative quadratic term. For smooth solutions the variational derivative \(w\) gives \(\mathrm d E/\mathrm dt=-\langle w,Lw\rangle\le0\).

Let \(X_h\) be the scalar cell-centered space and \(Y_h=X_h^m=(\mathbb R^m)^{I_{\mathrm{cell}}}\) the full auxiliary space, where
\[
 I_{\mathrm{cell}}=\prod_j\{1,\ldots,N_j\},\qquad
 h_j=\ell_j/N_j,\qquad x_{\boldsymbol i,j}=(i_j-\tfrac12)h_j .
\]
On either space, \(\langle v,w\rangle_h=(\prod_jh_j)\sum_{\boldsymbol i}v_{\boldsymbol i}\cdot w_{\boldsymbol i}\) and \(\|v\|_h^2=\langle v,v\rangle_h\). We use the same weight on each directional face grid. Periodic indexing gives the cell-to-face difference and its adjoint \cite{ShenEtAl2012,DongEtAl2018},
\[
 (D_{h,j}u)_{\boldsymbol i+e_j/2}
    =\frac{u_{\boldsymbol i+e_j}-u_{\boldsymbol i}}{h_j},
 \qquad
 (D_{h,j}^*z)_{\boldsymbol i}
    =\frac{z_{\boldsymbol i-e_j/2}-z_{\boldsymbol i+e_j/2}}{h_j}.
\]
The centered negative Laplacian is \(K_h=-\Delta_h=\sum_jD_{h,j}^*D_{h,j}\), and \(I_h\) denotes the discrete identity.

Let \(A_h,C_h,L_h\) be self-adjoint positive semidefinite operators on \(X_h\), and \(B_h:X_h\to Y_h\) linear. For \(D\) positive semidefinite write \(\|v\|_{D,h}^2=\langle Dv,v\rangle_h\). Define \(F_h(v_h)=\langle F(v_h),1\rangle_h\), so its gradient applies \(f\) at each cell. The discrete split and its flow are
\begin{gather}
 E_h(u_h)=\tfrac12\|u_h\|_{A_h,h}^2+F_h(B_hu_h)
                         -\tfrac12\|u_h\|_{C_h,h}^2,\label{eq:discrete-energy}\\
 \dot u_h=-L_hw_h,\qquad
 w_h=A_hu_h+B_h^*f(B_hu_h)-C_hu_h .
\end{gather}
Set \(T_h=\operatorname{ran}L_h=(\ker L_h)^\perp\) and let \(P_h\) be the orthogonal projection onto \(T_h\). Since \(\dot u_h\in T_h\), the flow preserves \(U_h=u_h^0+T_h\). The operator \(L_h^\dagger\) inverts \(L_h\) on \(T_h\) and vanishes on \(\ker L_h\).

\subsection{Convex splitting and time-step minimization}
\label{sec:time-step-solvers}
Let \(\tau>0\) and \(t_n=n\tau\). Backward Euler (BE) treats the convex part implicitly and the concave part explicitly \cite{Eyre1998}. This gives
\begin{equation}\label{eq:be-convex-splitting}
 \tfrac{u_h^{n+1}-u_h^n}{\tau}=-L_hw_h^{n+1},
 \qquad w_h^{n+1}=A_hu_h^{n+1}+B_h^*f(B_hu_h^{n+1})-C_hu_h^n.
\end{equation}
After one BE step we use the second-order backward differentiation formula (BDF2) with extrapolation of the explicit term and the following stabilization:
\begin{equation}\label{eq:bdf2-convex-splitting}
\begin{aligned}
 \tfrac{3u_h^{n+1}-4u_h^n+u_h^{n-1}}{2\tau}&=-L_hw_h^{n+1},\\
 w_h^{n+1}&=A_hu_h^{n+1}+B_h^*f(B_hu_h^{n+1})
           -C_h(2u_h^n-u_h^{n-1})+S_h(u_h^{n+1}-u_h^n),\\
 S_h&=\sigma\tau C_hL_hC_h,\qquad \sigma\ge0 .
\end{aligned}
\end{equation}
The operator \(S_h\) generalizes the Douglas--Dupont stabilization used for SS by Feng et~al.~\cite{FengEtAl2018}. In the SS case, \(L_h=I_h\) and \(C_h=K_h\) give \(S_h=\sigma\tau K_h^2\) in the chemical potential. The projection onto \(\ker L_h\) obeys the homogeneous BE or BDF2 recurrence, so the BE startup and the subsequent steps preserve the common initial kernel component. The stabilization is \(O(\tau^2)\) at fixed \(h\) on smooth trajectories and positive semidefinite for \(\sigma\ge0\). Applying \(L_h^\dagger\) gives the respective quadratic operators and loads
\begin{align}
 Q_{1,h}&=\tau^{-1}L_h^\dagger+A_h,
 & b_{1,h}^n&=(\tau^{-1}L_h^\dagger+C_h)u_h^n, \label{eq:be-step-coefficients}\\
 Q_{2,h}&=\tfrac{3}{2\tau}L_h^\dagger+A_h+S_h,
 & b_{2,h}^n&=\tau^{-1}L_h^\dagger(2u_h^n-\tfrac12u_h^{n-1})
             +C_h(2u_h^n-u_h^{n-1})+S_hu_h^n . \label{eq:bdf2-step-coefficients}
\end{align}
With \((Q_h,b_h^n)\) chosen accordingly, the step minimizes
\begin{equation}\label{eq:time-step-minimization}
 J_h(u_h)=g_h(u_h)+F_h(B_hu_h),\qquad
 g_h(u_h)=\tfrac12\langle Q_hu_h,u_h\rangle_h-\langle b_h^n,u_h\rangle_h,
 \qquad u_h\in U_h .
\end{equation}
Indeed, \(L_h^\dagger L_h=P_h\) gives the projected stationarity equation for \(J_h\). Write \(Q_{T,h}=P_hQ_h|_{T_h}\) and \(B_{T,h}=B_h|_{T_h}\), so \(B_{T,h}^*=P_hB_h^*\). Positivity of \(L_h^\dagger\) on \(T_h\) gives
\begin{equation}\label{eq:tangent-coercivity}
 \langle Q_{T,h}v_h,v_h\rangle_h\ge\mu_h\|v_h\|_h^2
 \quad(v_h\in T_h),\qquad \mu_h>0,
\end{equation}
and hence strong convexity on \(U_h\). Using parenthesized superscripts for inner iterates, define
\begin{equation}\label{eq:normalized-optimality-residual}
 r_h^{(k)}=P_h\bigl(Q_hu_h^{(k)}+B_h^*f(B_hu_h^{(k)})-b_h^n\bigr),
 \qquad R_h^{(k)}=\frac{\|r_h^{(k)}\|_h}{\max\{1,\|r_h^{(0)}\|_h\}} .
\end{equation}
The numerator is the norm of the projected gradient in the unpreconditioned grid inner product. The grid function \(L_hr_h^{(k)}\) is the corresponding defect in the time equation, since \(L_hP_h=L_h\).

\paragraph{Energy estimates}
For exact unforced steps, write \(a_h=u_h^{n+1}-u_h^n\), \(b_h=u_h^n-u_h^{n-1}\), and let \(\Delta E_h=E_h(u_h^{n+1})-E_h(u_h^n)\). Convexity bounds the increment of the implicit energy by its gradient at \(u_h^{n+1}\) paired with \(a_h\), and the explicit quadratic increment is \(\langle C_hu_h^n,a_h\rangle_h+\tfrac12\|a_h\|_{C_h,h}^2\). In BE this gives \(\Delta E_h+\tfrac12\|a_h\|_{C_h,h}^2\le\langle w_h^{n+1},a_h\rangle_h\). Testing the BE equation with \(L_h^\dagger a_h\) yields
\begin{equation}\label{eq:be-energy-dissipation}
 E_h(u_h^{n+1})-E_h(u_h^n)+\tfrac12\|a_h\|_{C_h,h}^2
      +\tau^{-1}\|a_h\|_{L_h^\dagger,h}^2\le0 .
\end{equation}
For BDF2, testing the discrete derivative \(\tau^{-1}[a_h+(a_h-b_h)/2]\) against \(L_h^\dagger a_h\) and using the same convexity estimate gives
\[
 \Delta E_h+\tfrac12\|a_h\|_{C_h,h}^2-\langle C_hb_h,a_h\rangle_h
 +\|a_h\|_{S_h,h}^2
 +\tau^{-1}\|a_h\|_{L_h^\dagger,h}^2
 +\tfrac{1}{2\tau}\langle L_h^\dagger(a_h-b_h),a_h\rangle_h\le0.
\]
Apply the polarization identity
\[
 2\langle D(a_h-b_h),a_h\rangle_h
 =\|a_h\|_{D,h}^2-\|b_h\|_{D,h}^2+\|a_h-b_h\|_{D,h}^2
\]
to the \(L_h^\dagger\) term and use \(\tfrac12\|a_h\|_{C_h,h}^2-\langle C_hb_h,a_h\rangle_h=\tfrac12\|a_h-b_h\|_{C_h,h}^2-\tfrac12\|b_h\|_{C_h,h}^2\). With
\begin{equation}\label{eq:modified-bdf2-energy}
 \widetilde E_h^n=E_h(u_h^n)+\tfrac{1}{4\tau}\|b_h\|_{L_h^\dagger,h}^2
                        +\tfrac12\|b_h\|_{C_h,h}^2 ,
\end{equation}
the exact unforced step satisfies \(\widetilde E_h^{n+1}-\widetilde E_h^n+D_{\mathrm{step}}\le0\), where
\begin{equation}\label{eq:bdf2-energy-estimate}
\begin{aligned}
 D_{\mathrm{step}}={}&\tfrac{1}{4\tau}\|a_h-b_h\|_{L_h^\dagger,h}^2
            +\tfrac12\|a_h-b_h\|_{C_h,h}^2+\tau^{-1}\|a_h\|_{L_h^\dagger,h}^2\\
 &+\sigma\tau\|L_h^{1/2}C_ha_h\|_h^2-\tfrac12\|a_h\|_{C_h,h}^2 .
\end{aligned}
\end{equation}
Since \(L_h^{1/2}(L_h^\dagger)^{1/2}=P_h\) and \(P_ha_h=a_h\), we have
\[
   \|a_h\|_{C_h,h}^2
   =\big\langle L_h^{1/2}C_ha_h,(L_h^\dagger)^{1/2}a_h\big\rangle_h .
\]
This identity uses the self-adjointness of the square roots and the projection identity. Young's inequality gives
\[
 \tfrac12\|a_h\|_{C_h,h}^2
 \le\tau^{-1}\|a_h\|_{L_h^\dagger,h}^2
     +\tfrac{\tau}{16}\|L_h^{1/2}C_ha_h\|_h^2 ,
\]
so \(D_{\mathrm{step}}\ge0\) for \(\sigma\ge1/16\). An inexact step adds the signed work \(W_{\mathrm{res}}=\langle r_h^{n+1},a_h\rangle_h\) on the right, using its unnormalized projected residual. In that estimate \(W_{\mathrm{res}}\le D_{\mathrm{step}}\) suffices for energy nonincrease.

\subsection{Model operators and symbols}
\label{sec:model-discretizations}
The slope-selection (SS) molecular beam epitaxy model for the film height \(u\), with \(\varepsilon_m>0\), is the \(L^2\) gradient flow of
\[
 E_{\mathrm{SS}}(u)=\tfrac{\varepsilon_m^2}{2}\|\Delta u\|_{L^2}^2
       +\tfrac14\int_\Omega(|\nabla u|^2-1)^2\,\mathrm dx .
\]
In \eqref{eq:continuous-gradient-flow}, \(L=I\) and \(w=\varepsilon_m^2\Delta^2u-\nabla\cdot[(|\nabla u|^2-1)\nabla u]\) \cite{LiLiu2003,FengEtAl2018}. The corresponding grid split uses \(L_h=I_h\), \(A_h=\varepsilon_m^2K_h^2\), \(C_h=K_h\), \(F(y)=|y|^4/4\), and a \(2d\)-component lifted slope
\begin{equation}\label{eq:mbe-lifted-gradient}
 (B_hu_h)_{\boldsymbol i}
 =2^{-1/2}\bigl((D_{h,j}u_h)_{\boldsymbol i+e_j/2},
                (D_{h,j}u_h)_{\boldsymbol i-e_j/2}\bigr)_{j=1}^d .
\end{equation}
Each face occurs twice in the cell sum with weight \(1/2\), giving \(\langle B_hu_h,B_hv_h\rangle_h=\langle K_hu_h,v_h\rangle_h\). Hence \(B_h^*B_h=K_h\). Here \(U_h=T_h=X_h\) and \(P_h=I_h\). Because \(A_h,C_h,S_h,B_h\) annihilate constants, \(\langle w_h^{n+1},1\rangle_h=0\), and exact BE startup followed by BDF2 conserves the sampled mean. Expanding the physical energy gives
\begin{equation}\label{eq:mbe-ss-discrete-energy}
 E_{h,\mathrm{SS}}^{\mathrm{phys}}(u_h)=\tfrac{\varepsilon_m^2}{2}\|K_hu_h\|_h^2
 +\tfrac14\langle (|B_hu_h|^2-1)^2,1\rangle_h=E_h(u_h)+|\Omega|/4 .
\end{equation}

The no-slope-selection (NSS) model is the \(L^2\) gradient flow of
\[
 E_{\mathrm{NSS}}(u)=\tfrac{\varepsilon_m^2}{2}\|\Delta u\|_{L^2}^2
       -\tfrac12\int_\Omega\log(1+|\nabla u|^2)\,\mathrm dx .
\]
It shares \(L_h\), \(A_h\), and the lifted gradient \(B_h\) with SS. Its split uses \(C_h=\beta K_h\) and \(F(y)=(\beta|y|^2-\log(1+|y|^2))/2\). The splitting \(E=E_c-E_e\) of Wang et~al.~\cite[Section~3.1]{WangWangWise2010} corresponds to \(\beta=1\). We take \(\beta>1\) by adding \(\int_\Omega(\beta-1)|\nabla u|^2/2\,\mathrm dx\) to both convex energies. These additions cancel in the physical energy, and the BDF2 stabilization depends on the modified \(C_h\).

NSS fits the general time discretization and ADMM framework. The mesh-uniform results in Theorem~\ref{thm:mesh-uniform-admm} and Corollary~\ref{cor:mesh-uniform-outer-trajectories} concern SS and standard PFC. The NSS curvature and penalty estimates are in \suppsection{sec:supp-penalty}{Penalty for globally bounded curvature}.

For the standard phase-field crystal (PFC) model, \(u\) is a conserved density and the \(H^{-1}\) gradient flow has unit preferred wavenumber \cite{ElderEtAl2002,WiseWangLowengrub2009}. Its chemical potential and evolution are
\[
 w=(I+\Delta)^2u+u^3-\epsilon u,\qquad \partial_tu=\Delta w ,
 \qquad \epsilon>0.
\]
The squared linear term and the quartic contribution form the convex part \cite[Eq.~(11)]{ShinLeeLee2016}. The local potential split also appears in the mixed formulation of Vignal et~al.~\cite[Section~3.1.3]{VignalEtAl2015}. The discrete operators are \(Y_h=X_h\), \(L_h=K_h\), \(A_h=(I_h-K_h)^2\), \(C_h=\epsilon I_h\), \(B_h=I_h\), and \(F(s)=s^4/4\). The discrete physical energy is
\begin{equation}\label{eq:pfc-discrete-energy}
 E_h(u_h)=\tfrac12\|(I_h-K_h)u_h\|_h^2+\tfrac14\langle u_h^4,1\rangle_h
                           -\tfrac{\epsilon}{2}\|u_h\|_h^2 .
\end{equation}
Let \(m_h=|\Omega|^{-1}\langle u_h^0,1\rangle_h\). The PFC affine space and tangent projection are
\[
 U_h=\{u_h\in X_h:\langle u_h,1\rangle_h=|\Omega|m_h\},\qquad
 P_hv_h=v_h-|\Omega|^{-1}\langle v_h,1\rangle_h\,1 .
\]
The zero mode of \(u_h\) is prescribed by this constraint. The inverse \(K_h^\dagger\) acts only on nonzero modes, and the auxiliary fields belong to the full space \(Y_h\). A feasible inexact PFC iterate retains the prescribed mean independently of its residual. Its time-equation defect is \(K_hr_h^{(k)}\).

On the periodic grid the symbol of \(K_h\) is
\begin{equation}\label{eq:mbe-k-symbol}
 \nu_{\boldsymbol k}=4\sum_{j=1}^d h_j^{-2}\sin^2(\pi k_j/N_j).
\end{equation}
For SS and PFC, writing \(\nu=\nu_{\boldsymbol k}\), \eqref{eq:be-step-coefficients}--\eqref{eq:bdf2-step-coefficients} give
\begin{equation}\label{eq:model-quadratic-symbols}
\begin{array}{c|cc}
 &\widehat Q_{1,h}&\widehat Q_{2,h}\\ \hline
 \text{SS}&\tau^{-1}+\varepsilon_m^2\nu^2
      &\tfrac{3}{2\tau}+(\varepsilon_m^2+\sigma\tau)\nu^2\\
 \text{PFC }(\nu>0)&(\tau\nu)^{-1}+(1-\nu)^2
      &\tfrac{3}{2\tau\nu}+(1-\nu)^2+\sigma\tau\epsilon^2\nu .
\end{array}
\end{equation}
The load multipliers in \eqref{eq:be-step-coefficients} and \eqref{eq:bdf2-step-coefficients} are
\[
\begin{array}{c|ccc}
 & \text{Time }(\widehat{L_h^\dagger}/\tau) & \text{Explicit }(\widehat C_h)
   & \text{Stabilization }(\widehat S_h)\\ \hline
 \text{SS}            & \tau^{-1}       & \nu      & \sigma\tau\nu^2\\
 \text{PFC }(\nu>0)   & (\tau\nu)^{-1}  & \epsilon & \sigma\tau\epsilon^2\nu .
\end{array}
\]
The first two columns act on \(u_h^n\) in BE. In BDF2 they act on \(2u_h^n-u_h^{n-1}/2\) and \(2u_h^n-u_h^{n-1}\), respectively, and the third column acts on \(u_h^n\).

The quartic Allen--Cahn and Cahn--Hilliard equations also fit \eqref{eq:time-step-minimization} with \(A_h=\varepsilon^2K_h\), \(B_h=C_h=I_h\), \(F(s)=s^4/4\), and \(L_h=I_h\) and \(L_h=K_h\), respectively. The latter uses fixed mean. ADMM methods for their Flory--Huggins variants treat singular endpoints and, for Cahn--Hilliard, use a convex--concave saddle-point splitting in the phase and chemical potential \cite{JiangLiangLu2027,LiLiangQiao2026}.

\section{ADMM formulation and fixed-point error bound}
\label{sec:theoretical-analysis}

\subsection{Predictor and accepted updates}

Split the nonlinear argument by imposing \(v_h=B_hu_h\) on \(g_h(u_h)+F_h(v_h)\), as in gradient splitting \cite{GlowinskiMarroco1975}. The augmented Lagrangian on \(U_h\times Y_h\) is
\begin{equation}\label{eq:augmented-step-lagrangian}
    \begin{aligned}
    \Lambda_{\rho,h}(u_h,v_h,\lambda_h)
        &=g_h(u_h)+F_h(v_h)+\langle\lambda_h,B_hu_h-v_h\rangle_h\\
        &\quad+\tfrac{\rho}{2}\|B_hu_h-v_h\|_h^2 .
    \end{aligned}
\end{equation}
Since \(J_h\) is coercive and strongly convex on \(U_h\), it has a unique minimizer \(u_h^*\). Its projected optimality condition gives the unique Karush--Kuhn--Tucker triple
\begin{equation}\label{eq:step-kkt}
    P_h(Q_hu_h^*-b_h^n+B_h^*\lambda_h^*)=0,\qquad
    x_h^*=v_h^*=B_hu_h^*,\qquad \lambda_h^*=f(v_h^*).
\end{equation}
Conversely, eliminating \(v_h^*,\lambda_h^*\) from \eqref{eq:step-kkt} recovers the optimality condition.

ADMM fixes the multiplier, minimizes \eqref{eq:augmented-step-lagrangian} in \(v_h\) and then \(u_h\), and updates the multiplier \cite{GabayMercier1976,BoydEtAl2011}. For the SS and PFC discretizations in Section~\ref{sec:model-discretizations}, the quadratic subproblem is solved by FFT on the admissible Fourier modes.

\begin{algorithm}[H]
\caption{ADMM for the time-step minimization problem}\label{alg:standard-admm}
\begin{subequations}\label{eq:admm-sweep}
\begin{enumerate}
    \item \label{step:admm-initialize} Choose \(\rho>0\), \(u_h^{(0)}\in U_h\), \(\lambda_h^{(0)}\in Y_h\), a tolerance \(\mathrm{tol}>0\), and an iteration limit \(K_{\max}\). Set \(k=0\).
    \item \label{step:admm-stop} Evaluate \(R_h^{(k)}\) from \eqref{eq:normalized-optimality-residual}. Return \(u_h^{(k)}\) if \(R_h^{(k)}\le\mathrm{tol}\). Otherwise, if \(k=K_{\max}\), report that the limit was reached and stop.
    \item \label{step:admm-proximal} Compute \(v_h^{(k+1)}\in Y_h\) cellwise from
          \begin{equation}\label{eq:admm-proximal-step}
              f(v_h^{(k+1)})+\rho(v_h^{(k+1)}-B_hu_h^{(k)})=\lambda_h^{(k)}.
          \end{equation}
    \item Find \(\bar u_h^{(k+1)}\in U_h\) by solving
          \begin{equation}\label{eq:admm-linear-step}
              P_h(Q_h+\rho B_h^*B_h)\bar u_h^{(k+1)}
              =P_h[b_h^n+B_h^*(\rho v_h^{(k+1)}-\lambda_h^{(k)})].
          \end{equation}
    \item \label{step:admm-multiplier} Form the multiplier predictor
          \begin{equation}\label{eq:admm-multiplier-step}
              \bar\lambda_h^{(k+1)}
              =\lambda_h^{(k)}+\rho(B_h\bar u_h^{(k+1)}-v_h^{(k+1)}).
          \end{equation}
    \item Accept the predictor with \(\alpha_k=1\), so \(u_h^{(k+1)}=\bar u_h^{(k+1)}\) and \(\lambda_h^{(k+1)}=\bar\lambda_h^{(k+1)}\). Increase \(k\) by one and return to Step~\ref{step:admm-stop}.
\end{enumerate}
\end{subequations}
\end{algorithm}

Earlier searches test a fixed-point residual \cite{GiselssonFaltBoyd2016} or the Douglas--Rachford envelope \cite{ThemelisStellaPatrinos2022}. ADMM-Search uses the original time-step objective along the primal predictor direction. See \suppsection{sec:supp-search}{Objective-based search} for the line-search implementation.

\begin{algorithm}[H]
\caption{ADMM-Search for the time-step minimization problem}\label{alg:admm-search}
\begin{enumerate}
    \item Initialize as in Step~\ref{step:admm-initialize} of Algorithm~\ref{alg:standard-admm}. Choose a fixed closed interval \([\alpha_{\min},\alpha_{\max}]\subset(0,2)\) containing \(1\).
    \item \label{step:admm-search-stop} Apply Step~\ref{step:admm-stop} of Algorithm~\ref{alg:standard-admm}.
    \item Compute the common predictor by Steps~\ref{step:admm-proximal}--\ref{step:admm-multiplier} of Algorithm~\ref{alg:standard-admm}. Form the direction
          \[
              d_h^{(k)}=\bar u_h^{(k+1)}-u_h^{(k)}.
          \]
    \item Set \(\alpha_k=1\), giving an unrelaxed first sweep (\(\alpha_0=1\)). For \(k\ge1\) and \(d_h^{(k)}\ne0\), search the interval for a candidate \(\alpha\) with \(\phi_k(\alpha)<\phi_k(1)\), where \(\phi_k(\alpha)=J_h(u_h^{(k)}+\alpha d_h^{(k)})\). Set \(\alpha_k=\alpha\) only if the search finds such a candidate. Retain the unit step if the direction is zero, the search fails, or the candidate gives no strict decrease.
    \item Retain \(v_h^{(k+1)}\) and update the primal variable and multiplier with the same \(\alpha_k\):
          \begin{equation}\label{eq:admm-search-update}
              (u_h,\lambda_h)^{(k+1)}=(1-\alpha_k)(u_h,\lambda_h)^{(k)}
              +\alpha_k(\bar u_h,\bar\lambda_h)^{(k+1)}.
          \end{equation}
          Increase \(k\) by one and return to Step~\ref{step:admm-search-stop}.
\end{enumerate}
\end{algorithm}

\subsection{Relaxed Douglas--Rachford iteration}

For the analysis we fix \(h\), the time level and the penalty \(\rho>0\). We continue the iteration in exact arithmetic with exact subproblem solves and no tolerance or iteration cap. The reduction below uses the classical ADMM--Douglas--Rachford (DR) equivalence \cite{LionsMercier1979,EcksteinBertsekas1992}. Write \(s_h(u,\lambda)=P_h(Q_hu-b_h^n+B_h^*\lambda)\) for the projected quadratic stationarity map.

\begin{lemma}[Reduction to a fixed-point iteration]\label{lem:admm-dr-iteration}
Let \(x_h^{(k)}=B_hu_h^{(k)}\) and \(z_h^{(k)}=x_h^{(k)}-\lambda_h^{(k)}/\rho\). For \(k\geq1\), the accepted iterates satisfy
\begin{equation}\label{eq:admm-affine-stationarity}
    P_h(Q_hu_h^{(k)}-b_h^n+B_h^*\lambda_h^{(k)})=0.
\end{equation}
For \(x\in B_hU_h\), set \(G_h(x)=\min\{g_h(u):u\in U_h,\ B_hu=x\}\) and extend by \(+\infty\) elsewhere in \(Y_h\). Write \(P_H=\operatorname{prox}_{H/\rho}\) and \(R_H=2P_H-I_{Y_h}\). The state satisfies the relaxed Douglas--Rachford (DR) iteration
\begin{equation}\label{eq:admm-relaxed-dr}
    \begin{aligned}
        x_h^{(k)}&=P_{G_h}z_h^{(k)},\qquad
        v_h^{(k+1)}=P_{F_h}R_{G_h}z_h^{(k)},\\
        z_h^{(k+1)}&=z_h^{(k)}+\alpha_k(v_h^{(k+1)}-x_h^{(k)})
          =\bigl[(1-\alpha_k/2)I_{Y_h}+(\alpha_k/2)R_{F_h}R_{G_h}\bigr]z_h^{(k)}.
    \end{aligned}
\end{equation}
\end{lemma}

\begin{proof}
Combining \eqref{eq:admm-linear-step} with \eqref{eq:admm-multiplier-step} gives \(s_h(\bar u_h^{(k+1)},\bar\lambda_h^{(k+1)})=0\). The accepted pair satisfies \eqref{eq:admm-search-update}, with \(\alpha_k=1\) for standard ADMM. Since \(s_h\) is affine,
\[
    s_h(u_h^{(k+1)},\lambda_h^{(k+1)})
      =(1-\alpha_k)s_h(u_h^{(k)},\lambda_h^{(k)}).
\]
The first sweep establishes \eqref{eq:admm-affine-stationarity}, and the common update preserves it. Substitute \(\lambda_h^{(k)}=\rho(B_hu_h^{(k)}-z_h^{(k)})\) to obtain the constrained optimality condition for the unique recovery
\begin{equation}\label{eq:admm-primal-recovery}
    u_h^{(k)}=\operatorname*{argmin}_{u\in U_h}
      \bigl\{g_h(u)+(\rho/2)\|B_hu-z_h^{(k)}\|_h^2\bigr\}.
\end{equation}
The fibers in the definition of \(G_h\) have unique quadratic minimizers. Its sublevel sets are images of compact sublevel sets of \(g_h\) on \(U_h\), so \(G_h\) is proper, closed and convex. Minimizing \eqref{eq:admm-primal-recovery} first over each fiber reduces it to the minimization of \(G_h(x)+(\rho/2)\|x-z_h^{(k)}\|_h^2\), which gives the first proximal map. The optimality condition \eqref{eq:admm-proximal-step} gives the second proximal map at \(x_h^{(k)}+\lambda_h^{(k)}/\rho=2x_h^{(k)}-z_h^{(k)}\). Finally, \eqref{eq:admm-multiplier-step} gives
\[
    B_h\bar u_h^{(k+1)}-\bar\lambda_h^{(k+1)}/\rho
      =z_h^{(k)}+v_h^{(k+1)}-x_h^{(k)} .
\]
Relaxing the pair in this identity yields \eqref{eq:admm-relaxed-dr}.
\end{proof}

At any fixed point \(z_h\) of \(T=R_{F_h}R_{G_h}\), the two proximal values coincide because \(Tz_h-z_h=2(v_h-x_h)\). Recover \(u_h\) by \eqref{eq:admm-primal-recovery} and set \(\lambda_h=\rho(x_h-z_h)\). The recovery condition gives \(s_h(u_h,\lambda_h)=0\), and the second proximal condition gives \(\lambda_h=f(v_h)\). Thus the recovered triple solves \eqref{eq:step-kkt}. Conversely, these optimality conditions yield the two proximal relations at \(z_h^*=x_h^*-\lambda_h^*/\rho\), with \(x_h^*=v_h^*\). Uniqueness of the KKT triple gives uniqueness of the fixed point.

The reflections are nonexpansive by firm nonexpansiveness of proximal maps \cite{Moreau1965}. With \(\theta=\alpha_k/2\), the relaxed state is a convex combination of \(z_h^{(k)}\) and \(Tz_h^{(k)}\). Expanding the squared error and using \(\|Tz_h^{(k)}-Tz_h^*\|_h\le\|z_h^{(k)}-z_h^*\|_h\) leaves the decrement \(\theta(1-\theta)\|Tz_h^{(k)}-z_h^{(k)}\|_h^2\). Hence
\begin{equation}\label{eq:admm-fejer}
    \|z_h^{(k+1)}-z_h^*\|_h^2
       \leq \|z_h^{(k)}-z_h^*\|_h^2
          -\alpha_k(2-\alpha_k)\|x_h^{(k)}-v_h^{(k+1)}\|_h^2 .
\end{equation}

\subsection{Error bound and fixed-grid convergence}

On a fixed grid the Fej\'er bound makes the state orbit bounded. The proximal maps and reflections in \eqref{eq:admm-relaxed-dr} then bound the two proximal orbits. If \(f\) is locally Lipschitz, it has a Lipschitz constant on a compact convex set containing both orbits and the solution. The following estimate identifies the relevant dependence on the quadratic block.

\begin{lemma}[Fixed-point error bound]\label{lem:admm-error-bound}
Set \(\kappa_h=\|B_{T,h}Q_{T,h}^{-1}B_{T,h}^*\|\). If \(f\) has Lipschitz constant \(L_{\rm loc}\) on a convex set containing \(x_h^{(k)},v_h^{(k+1)},x_h^*\), then
\begin{equation}\label{eq:admm-error-bound}
    \|z_h^{(k)}-z_h^*\|_h\le C_{\rho,h}\|x_h^{(k)}-v_h^{(k+1)}\|_h,\quad
    C_{\rho,h}=(1+L_{\rm loc}/\rho)\bigl[1+\kappa_h(L_{\rm loc}+\rho)\bigr].
\end{equation}
\end{lemma}

\begin{proof}
Suppress the subscript \(h\) on the iterates. Set \(\xi=u^{(k)}-u^*\), \(x=x^{(k)}\), \(v=v^{(k+1)}\), and \(e=x-v\). The proximal condition gives \(\lambda^{(k)}=f(v)-\rho e\). Subtracting \eqref{eq:step-kkt} from \eqref{eq:admm-affine-stationarity} and testing by \(\xi\) gives
\[
    \begin{aligned}
    \|\xi\|_{Q_{T,h},h}^2
      &=-\langle f(x)-f(x^*),x-x^*\rangle_h
         +\langle f(x)-f(v)+\rho e,x-x^*\rangle_h\\
      &\le (L_{\rm loc}+\rho)\|e\|_h\|x-x^*\|_h,
    \end{aligned}
\]
where \(\|\xi\|_{Q_{T,h},h}^2=\langle Q_{T,h}\xi,\xi\rangle_h\) and the first term is nonpositive by monotonicity. From \(x-x^*=B_{T,h}\xi\) and \(\|B_{T,h}Q_{T,h}^{-1/2}\|^2=\kappa_h\), we have \(\|x-x^*\|_h\le\sqrt{\kappa_h}\|\xi\|_{Q_{T,h},h}\). Divide the resulting bound on \(\|\xi\|_{Q_{T,h},h}^2\) by \(\|\xi\|_{Q_{T,h},h}\) if \(\xi\ne0\), and apply the same inequality once more. This gives \(\|x-x^*\|_h\le\kappa_h(L_{\rm loc}+\rho)\|e\|_h\), which also holds if \(\xi=0\). Finally, \(z^{(k)}-z^*=x-x^*+e-(f(v)-f(x^*))/\rho\) and \(\|v-x^*\|_h\le\|x-x^*\|_h+\|e\|_h\) yield
\[
    \|z^{(k)}-z^*\|_h\le (1+L_{\rm loc}/\rho)(\|x-x^*\|_h+\|e\|_h).
\]
Inserting the bound on \(\|x-x^*\|_h\) proves \eqref{eq:admm-error-bound}.
\end{proof}

Combining a fixed-point error bound with Fej\'er decrease is a standard route to linear convergence \cite{LiuYuanZengZhang2018,PenaVeraZuluaga2021}. Here the preceding argument supplies the bound in terms of \(\kappa_h\) and \(L_{\rm loc}\).

\begin{theorem}[Convergence on a fixed grid]\label{thm:admm-linear}
With an unrelaxed first sweep (\(\alpha_0=1\)) and subsequent synchronized parameters in the prescribed interval, the iterates converge to \eqref{eq:step-kkt} for each \(\rho>0\). If \(f\) is locally Lipschitz, the state converges Q-linearly with factor at most
\begin{equation}\label{eq:admm-local-contraction}
    (1-\eta/C_{\rho,h}^2)^{1/2},\qquad
    \eta=\min_{\alpha\in[\alpha_{\min},\alpha_{\max}]}\alpha(2-\alpha)>0.
\end{equation}
The recovered variables and the stopping residual have geometric envelopes with the same factor. For standard ADMM, \(\eta=1\).
\end{theorem}
\begin{proof}
Summing \eqref{eq:admm-fejer} shows that
\(\eta\sum_{k\ge1}\|x_h^{(k)}-v_h^{(k+1)}\|_h^2\le\|z_h^{(1)}-z_h^*\|_h^2\). The state is bounded and the defect tends to zero. By continuity of the proximal maps, every cluster point is the unique fixed point \(z_h^*\). Thus the full state sequence converges. With local Lipschitz continuity, \eqref{eq:admm-error-bound} bounds the squared defect below by \(C_{\rho,h}^{-2}\|z_h^{(k)}-z_h^*\|_h^2\). Using this in \eqref{eq:admm-fejer} gives \eqref{eq:admm-local-contraction}.

Subtracting the optimality equations for the recoveries from \(z_h,\widetilde z_h\) gives
\[
    (Q_{T,h}+\rho B_{T,h}^*B_{T,h})(u_h-\widetilde u_h)
        =\rho B_{T,h}^*(z_h-\widetilde z_h).
\]
Thus the recovery is affine with Lipschitz constant at most \(\rho\|B_{T,h}\|/\mu_h\) by \eqref{eq:tangent-coercivity}. The second proximal formula and \(\lambda_h^{(k)}=\rho(x_h^{(k)}-z_h^{(k)})\) give Lipschitz dependence of \(v_h^{(k+1)}\) and \(\lambda_h^{(k)}\) on the state, yielding convergence of the recovered iterates and their geometric envelopes. The projected gradient vanishes at \(u_h^*\) and is Lipschitz on the bounded primal orbit. Applying \eqref{eq:normalized-optimality-residual}, with its fixed denominator, transfers the envelope to the stopping residual.
\end{proof}

\section{Mesh-uniform convergence and exact trajectories}
\label{sec:mesh-uniform-convergence}

\subsection{Single-step estimate}
\label{sec:mesh-uniform-result}

On a fixed periodic box \(\Omega=\prod_{j=1}^d(0,\ell_j)\), \(d\le3\), write
\(\|w_h\|_{H_h^s}=\|(I_h+K_h)^{s/2}w_h\|_h\), with \(K_h\) acting componentwise on lifted fields. We denote the weighted grid \(L^p\) norm by \(\|\cdot\|_{p,h}\) and the maximum norm by \(\|\cdot\|_{\infty,h}\), using Euclidean norms for vector values. All estimates below allow independent widths \(h_j=\ell_j/N_j\).

\begin{theorem}[Mesh-uniform convergence for one time step]
\label{thm:mesh-uniform-admm}
Consider an unforced BE or stabilized BDF2 step for SS or standard PFC with \(\sigma\ge0\) and \(0<\tau\le\tau_0\). Fix the model parameters, \(\rho>0\), and \(0<\alpha_{\min}\le1\le\alpha_{\max}<2\). Assume the history fields are uniformly bounded in \(H_h^2\) for PFC and \(H_h^3\) for SS, with the prescribed mean in PFC. Solve the subproblems exactly, take an unrelaxed first sweep (\(\alpha_0=1\)), and use a common relaxation \(\alpha_k\in[\alpha_{\min},\alpha_{\max}]\) for \(k\ge1\). Suppose \(z_h^{(1)}\) is uniformly bounded in \(\|\cdot\|_{\infty,h}\) and, for SS, in \(H_h^1\). Then for all \(k\ge1\),
\begin{equation}\label{eq:mesh-uniform-residual-bound}
  \|z_h^{(k+1)}-z_h^\ast\|_h
     \le q_{\mathrm{DR}}\|z_h^{(k)}-z_h^\ast\|_h,\qquad
  \|r_h^{(k)}\|_h\le C_R q_{\mathrm{DR}}^{k-1},
\end{equation}
where \(0<q_{\mathrm{DR}}<1\) and \(C_R>0\) are independent of \(h\), \(k\), and \(0<\tau\le\tau_0\). The denominator in \eqref{eq:normalized-optimality-residual} is at least one, so \(R_h^{(k)}\le C_R q_{\mathrm{DR}}^{k-1}\).
\end{theorem}

Let \(\delta_{h,j}w_h=(w_h(\cdot+h_je_j)-w_h)/h_j\), and collect these differences in \(\delta_h w_h\). Parseval gives
\[
   \|w_h\|_{H_h^1}^2=\|w_h\|_h^2+\|\delta_hw_h\|_h^2,
\]
with the squared grid norm summing over directions and components on the right. Let \(\mathcal K_h\) contain the principal Fourier indices \(-\lfloor N_j/2\rfloor\le k_j<\lceil N_j/2\rceil\). On these modes, the symbol \(\nu_{\boldsymbol k}\) in \eqref{eq:mbe-k-symbol} satisfies
\[
  \frac4{\pi^2}|\boldsymbol\xi_{\boldsymbol k}|^2
     \le \nu_{\boldsymbol k}\le|\boldsymbol\xi_{\boldsymbol k}|^2,
  \qquad \xi_{\boldsymbol k,j}=2\pi k_j/\ell_j .
\]
This follows coordinate by coordinate from \(2|t|/\pi\le|\sin t|\le |t|\) for \(|t|\le\pi/2\). The spectral norm is therefore uniformly equivalent to the Fourier norm with weight \((1+|\boldsymbol\xi_{\boldsymbol k}|^2)^s\), independent of the ratios of the grid widths. The embeddings used below are
\begin{equation}\label{eq:mesh-discrete-embeddings}
  \|w_h\|_{\infty,h}\le C_\Omega\|w_h\|_{H_h^2},
  \qquad
  \|w_h\|_{6,h}\le C_\Omega\|w_h\|_{H_h^1}.
\end{equation}
For the first, Cauchy--Schwarz in the inverse Fourier series gives
\[
  \|w_h\|_{\infty,h}\le
   |\Omega|^{-1/2}
   \left(\sum_{\boldsymbol k\in\mathcal K_h}(1+\nu_{\boldsymbol k})^{-2}\right)^{1/2}
   \|w_h\|_{H_h^2}.
\]
The frequency comparison bounds the sum by a multiple of
\(\sum_{\boldsymbol k\in\mathbb Z^d}(1+|\boldsymbol k|^2)^{-2}<\infty\) for \(d\le3\).
See \cite[Lemma~4.1]{EvenDarMandelSchochet2017} for discrete Fourier analysis and \cite[Lemma~3.1, Eqs.~(3.5) and (B.23)]{GuoWangWiseYue2016} for the cell-centered estimate with equal widths.

For the second estimate, let \(w\) be the continuous periodic multilinear interpolant of the cell values. Reference-cube norm equivalence gives
\[
  \|w_h\|_{6,h}\le C\|w\|_{L^6(\Omega)},\qquad
  \|w\|_{L^2(\Omega)}\le C\|w_h\|_h,\qquad
  \|\partial_j w\|_{L^2(\Omega)}\le C\|\delta_{h,j}w_h\|_h .
\]
In the last bound, the derivative interpolates the physical edge differences in the other directions. On scaling each coordinate of a cell separately, the only measure factor is the cell volume. The factor \(1/h_j\) in the derivative is already in the edge difference. Thus the constants are independent of all widths and their ratios. The continuous Sobolev embedding \cite[Section~5.6]{Evans2010} gives the second bound in \eqref{eq:mesh-discrete-embeddings}. See \cite[Lemma~C.1, Eq.~(C.1)]{ChenWangWise2022} for the equal-width cell-centered bound.

\subsection{Uniform orbit and state contraction}
\label{sec:mesh-uniform-orbit}

\begin{proof}[State contraction in Theorem~\ref{thm:mesh-uniform-admm}]
For the quadratic symbols in \eqref{eq:model-quadratic-symbols}, we have
\begin{equation}\label{eq:uniform-spectral-lower-bound}
   \widehat Q_h(\nu)\ge c_0(1+\nu)^2,\qquad
   c_0=
    \begin{cases}
     \frac12\min\{\tau_0^{-1},\varepsilon_m^2\},&\text{SS},\\
     \frac19\min\{1,(2\tau_0)^{-1}\},&\text{PFC}.
    \end{cases}
\end{equation}
For SS this holds on all modes since \(1+\nu^2\ge(1+\nu)^2/2\). For PFC, the time term is at least \(1/(2\tau_0)\) for \(0<\nu\le2\). For \(\nu\ge2\), \((1-\nu)^2\ge(1+\nu)^2/9\). The PFC bound applies on the tangent space, with the mean treated by the affine constraint. For SS the modal norm of \(B_{T,h}Q_{T,h}^{-1}B_{T,h}^\ast\) is \(\nu/\widehat Q_h\), since \(B_h^\ast B_h=K_h\). For PFC it is \(1/\widehat Q_h\) on nonzero modes and zero on the mean. Thus \(\kappa_h\le\bar\kappa\), with \(\bar\kappa=1/(4c_0)\) for SS and \(1/c_0\) for PFC.

Let \(M_{\rm hist}\) bound the prescribed history norms, with \(s=3\) for SS and \(s=2\) for PFC. The quadratic minimizer \(u_{g,h}=\arg\min_{U_h}g_h\) has coefficients \(\widehat b_h^n/\widehat Q_h\) on tangent modes. For PFC its constant mode is prescribed and has grid norm at most \(M_{\rm hist}\). In \eqref{eq:be-step-coefficients} and \eqref{eq:bdf2-step-coefficients}, the time coefficient divided by \(\widehat Q_h\) is at most one. Its history factor is \(u_h^n\) for BE and \((4u_h^n-u_h^{n-1})/3\) for BDF2. The stabilization coefficient divided by \(\widehat Q_h\) is also at most one and multiplies \(u_h^n\). The explicit coefficient divided by \(\widehat Q_h\) has bound \(c_{\rm exp}=1/(4c_0)\) for SS and \(c_{\rm exp}=\epsilon/c_0\) for PFC, by \eqref{eq:uniform-spectral-lower-bound}. It multiplies \(u_h^n\) in BE and \(2u_h^n-u_h^{n-1}\) in BDF2. Applying these multiplier bounds in \(H_h^s\) yields
\[
  \|u_{g,h}\|_{H_h^s}\le
  \begin{cases}
      (2+c_{\rm exp})M_{\rm hist},&\text{BE},\\
      (11/3+3c_{\rm exp})M_{\rm hist},&\text{BDF2}.
  \end{cases}
\]
The bounds include the possible constant mode in PFC.
For the SS lift, the symbol of \(B_h\) is a column with norm \(\sqrt\nu\), and the adjoint symbol is its conjugate transpose. Hence
\[
  \|B_hw_h\|_{H_h^s}\le\|w_h\|_{H_h^{s+1}},\qquad
  \|B_h^\ast v_h\|_{H_h^{-1}}\le\|v_h\|_h .
\]
In both models, \(B_hu_{g,h}\) is bounded in \(H_h^2\), and thus in the maximum norm. The energy comparison \(J_h(u_h^\ast)\le J_h(u_{g,h})\), \(F_h\ge0\), and the quadratic expansion at \(u_{g,h}\) give
\[
  \tfrac{c_0}{2}\|u_h^\ast-u_{g,h}\|_{H_h^2}^2
  \le g_h(u_h^\ast)-g_h(u_{g,h})
  \le F_h(B_hu_{g,h})\le C .
\]
It follows that \(x_h^\ast=B_hu_h^\ast\) is bounded in \(H_h^1\) for SS and \(H_h^2\) for PFC. From the quartic gradient and \eqref{eq:mesh-discrete-embeddings},
\(\|\lambda_h^\ast\|_h=\|x_h^\ast\|_{6,h}^3\le C\). Consequently \(z_h^\ast=x_h^\ast-\lambda_h^\ast/\rho\) is uniformly bounded in the grid norm. Fej\'er monotonicity \eqref{eq:admm-fejer} and the first-sweep bound yield
\begin{equation}\label{eq:initial-state-bound}
  \|z_h^{(k)}-z_h^\ast\|_h
  \le\|z_h^{(1)}-z_h^\ast\|_h\le D_\star,\qquad k\ge1 ,
\end{equation}
with \(D_\star\) independent of \(h\) and \(\tau\).

For \((u,z)=(u_h^{(k)},z_h^{(k)})\) or \((u_h^\ast,z_h^\ast)\), subtract the stationarity equation for \(u_{g,h}\) from the recovery equation \eqref{eq:admm-primal-recovery} to obtain
\[
  (Q_{T,h}+\rho B_{T,h}^\ast B_{T,h})(u-u_{g,h})
        =\rho B_{T,h}^\ast(z-B_hu_{g,h}).
\]
Set \(M_h=Q_{T,h}+\rho B_{T,h}^\ast B_{T,h}\) on \(T_h\). The added term is nonnegative and diagonal in the scalar Fourier representation. At each tangent mode and for every \(s\in\mathbb R\),
\[
   (1+\nu)^{s+4}|\widehat{(M_h^{-1}w)}_{\boldsymbol k}|^2
      \le c_0^{-2}(1+\nu)^s|\widehat w_{\boldsymbol k}|^2 .
\]
Summing yields \(\|M_h^{-1}w\|_{H_h^{s+4}}\le c_0^{-1}\|w\|_{H_h^s}\). In PFC, \(B_{T,h}^\ast=P_h\), which is nonexpansive in the grid norm. In SS, \(B_{T,h}^\ast=B_h^\ast\) maps the grid space to \(H_h^{-1}\) with norm at most one. The recovery equation therefore yields
\[
 \begin{aligned}
 \|u-u_{g,h}\|_{H_h^4}&\le(\rho/c_0)\|z-B_hu_{g,h}\|_h
                                  &&\text{for PFC},\\
 \|u-u_{g,h}\|_{H_h^3}&\le(\rho/c_0)\|z-B_hu_{g,h}\|_h
                                  &&\text{for SS}.
 \end{aligned}
\]
The right sides are uniformly bounded by \eqref{eq:initial-state-bound} and the history bound on \(u_{g,h}\). Using the slope symbol estimate for SS gives \(H_h^2\) and pointwise bounds on \(x_h^{(k)}=B_hu_h^{(k)}\) and \(x_h^\ast\) in both models.

Let \(U\ge\sqrt\rho\) be a uniform bound for these pointwise norms. At a cell, write \(x,z,v,z^+\) for the values in one step of \eqref{eq:admm-relaxed-dr}, and put \(p_\rho=\operatorname{prox}_{F/\rho}\). The quartic proximal update satisfies
\[
 v=p_\rho(2x-z)=s(2x-z),\qquad s=(1+|v|^2/\rho)^{-1},\qquad
 z^+=(1-\alpha s)z+\alpha(2s-1)x .
\]
If \(|v|\ge U\), then \(s\le1/2\) and \(1-\alpha s\ge0\). Since \(|v|\le s(2U+|z|)\), we have
\[
  |z^+|\le(1-\alpha s)|z|+\alpha(1-2s)U
     \le |z|+\alpha(U-|v|)\le |z|.
\]
If \(|v|<U\), substitute \(z=2x-v-|v|^2v/\rho\) to find
\[
  z^+=(2-\alpha)x+(\alpha-1)v-|v|^2v/\rho,\qquad
  |z^+|\le3U+U^3/\rho.
\]
Let \(Z_1\) be the uniform bound on \(\|z_h^{(1)}\|_{\infty,h}\) and put \(Z=\max\{Z_1,3U+U^3/\rho\}\). Induction gives \(\|z_h^{(k)}\|_{\infty,h}\le Z\) for \(k\ge1\). The fixed point satisfies the same bound since \(z_h^\ast=x_h^\ast-f(x_h^\ast)/\rho\). Both inputs and outputs of \(p_\rho\) are therefore bounded by \(M=2U+Z\), using nonexpansiveness and \(p_\rho(0)=0\). The quartic Hessian obeys \(\|Df(y)\|\le3|y|^2\), and \(L_\star=3M^2\) is a common Lipschitz constant of \(f\) on the resulting pointwise ball, hence also on the corresponding grid set in the grid norm.

With the uniform constants
\begin{equation}\label{eq:uniform-dr-factor}
 \begin{aligned}
 C_{\mathrm{EB}}&=(1+L_\star/\rho)[1+\bar\kappa(L_\star+\rho)],\\
 \eta_\alpha&=\min_{\alpha\in[\alpha_{\min},\alpha_{\max}]}\alpha(2-\alpha),
 \qquad q_{\mathrm{DR}}=\sqrt{1-\eta_\alpha/C_{\mathrm{EB}}^2},
 \end{aligned}
\end{equation}
the error bound \eqref{eq:admm-error-bound} and the Fej\'er inequality \eqref{eq:admm-fejer} establish the asserted contraction.
\end{proof}

\subsection{Transfer to the projected residual}
\label{sec:mesh-uniform-stopping-residual}

\begin{proof}[Residual estimate in Theorem~\ref{thm:mesh-uniform-admm}]
The affine stationarity relation \eqref{eq:admm-affine-stationarity} gives, for \(k\ge1\),
\begin{equation}\label{eq:stationary-residual-identity}
 r_h^{(k)}=B_{T,h}^\ast\bigl(f(x_h^{(k)})-\lambda_h^{(k)}\bigr),
 \qquad f(x_h^\ast)=\lambda_h^\ast .
\end{equation}
Put \(e_k=z_h^{(k)}-z_h^\ast\) and \(d_k=x_h^{(k)}-x_h^\ast\). The state contraction gives \(\|e_k\|_h\le D_\star q_{\mathrm{DR}}^{k-1}\). For PFC, \(B_{T,h}^\ast=P_h\) has norm at most one. Both \(\operatorname{prox}_{G_h/\rho}\) and its complement are nonexpansive. Since \(\lambda_h^{(k)}=\rho(x_h^{(k)}-z_h^{(k)})\), these two bounds give
\[
  \|d_k\|_h\le\|e_k\|_h,\qquad
  \|\lambda_h^{(k)}-\lambda_h^\ast\|_h\le\rho\|e_k\|_h .
\]
Inserting these estimates in \eqref{eq:stationary-residual-identity} gives \(\|r_h^{(k)}\|_h\le(L_\star+\rho)D_\star q_{\mathrm{DR}}^{k-1}\).

For SS, subtract the two recovery equations in \eqref{eq:admm-primal-recovery}. The lifted difference is
\[
  d_k=\rho B_h(Q_h+\rho K_h)^{-1}B_h^\ast e_k .
\]
The symbol of this map has norm \(\rho\nu/(\widehat Q_h+\rho\nu)\), bounded by \((\rho/c_0)/(1+\nu)\) from \eqref{eq:uniform-spectral-lower-bound}. Hence recovery gains two orders from \(L_h^2\). Using \eqref{eq:mesh-discrete-embeddings} gives
\begin{equation}\label{eq:ss-recovery-smoothing}
  \|d_k\|_{H_h^2}+\|d_k\|_{\infty,h}\le Cq_{\mathrm{DR}}^{k-1}.
\end{equation}
We use this smoothing to control the first differences before applying the slope adjoint.

On the input ball, the symmetric matrix
\[
 Dp_\rho(y)=\bigl(I+\rho^{-1}Df(p_\rho(y))\bigr)^{-1}
\]
has eigenvalues in \([m_p,1]\), \(m_p=(1+L_\star/\rho)^{-1}\). For the cellwise update \(\Phi_\alpha(x,z)=z+\alpha[p_\rho(2x-z)-x]\), \(D_z\Phi_\alpha=I-\alpha Dp_\rho\) and \(D_x\Phi_\alpha=\alpha(2Dp_\rho-I)\). Consequently
\[
 \|D_z\Phi_\alpha\|\le\gamma_p,\qquad
 \|D_x\Phi_\alpha\|\le\alpha_{\max},\qquad
 \gamma_p=\max\{1-\alpha_{\min}m_p,\alpha_{\max}-1\}.
\]
The factor \(\gamma_p\) is strictly smaller than \(q_{\mathrm{DR}}\). In fact \(C_{\mathrm{EB}}\ge m_p^{-1}\) and \(0<m_p<1\), so
\[
 \begin{aligned}
  1-q_{\mathrm{DR}}^2
     &\le \alpha_{\min}(2-\alpha_{\min})m_p^2
       <\alpha_{\min}m_p(2-\alpha_{\min}m_p),\\
  1-q_{\mathrm{DR}}^2
     &\le \alpha_{\max}(2-\alpha_{\max})m_p^2
       <\alpha_{\max}(2-\alpha_{\max}).
 \end{aligned}
\]
These bounds compare the squares of the two terms in \(\gamma_p\) with \(q_{\mathrm{DR}}^2\).

Compare the accepted and fixed-point updates at each cell, and then compare adjacent cells at the same sweep. The common scalar \(\alpha_k\) gives
\begin{equation}\label{eq:ss-local-recurrences}
 \begin{aligned}
  \|e_{k+1}\|_{\infty,h}&\le\gamma_p\|e_k\|_{\infty,h}
                      +\alpha_{\max}\|d_k\|_{\infty,h},\\
  \|\delta_h z_h^{(k+1)}\|_h&\le\gamma_p\|\delta_h z_h^{(k)}\|_h
                      +\alpha_{\max}\|\delta_h x_h^{(k)}\|_h .
 \end{aligned}
\end{equation}
A recurrence \(a_{k+1}\le\gamma_p a_k+C q_{\mathrm{DR}}^{k-1}\) gives
\(a_k\le q_{\mathrm{DR}}^{k-1}(a_1+C/(q_{\mathrm{DR}}-\gamma_p))\) by geometric summation. The first line of \eqref{eq:ss-local-recurrences} and \eqref{eq:ss-recovery-smoothing} therefore give \(\|e_k\|_{\infty,h}\le Cq_{\mathrm{DR}}^{k-1}\). The second line and the uniform \(H_h^2\) bound on \(x_h^{(k)}\) give a uniform \(H_h^1\) bound on \(z_h^{(k)}\). The fixed point \(z_h^\ast=x_h^\ast-f(x_h^\ast)/\rho\) has the same regularity bound by taking first differences and using the Lipschitz bound of \(f\) on the range of \(x_h^\ast\).

Write \(y_k=2x_h^{(k)}-z_h^{(k)}\) and \(y_\ast=2x_h^\ast-z_h^\ast\). Set
\[
  A_k=\int_0^1 Dp_\rho((1-t)y_\ast+t y_k)\,dt
\]
cellwise. Its eigenvalues lie in \([m_p,1]\). For a matrix field \(A\), \(\|A\|_{h,\mathrm F}\) is the weighted grid \(L^2\) norm of its cellwise Frobenius norms, summed in squares over directions for \(\delta_h A\). The identity \(Df(v)=|v|^2I+2vv^\mathsf T\) gives \(\operatorname{Lip}(Df)\le(2\sqrt{2d}+4)M\) in Frobenius norm on the ball of radius \(M\). Here the two contributions use \(\bigl||v|^2-|w|^2\bigr|\le2M|v-w|\) and the expansion \(vv^\mathsf T-ww^\mathsf T=(v-w)v^\mathsf T+w(v-w)^\mathsf T\). The inverse difference identity, \(\|Dp_\rho\|\le1\), and nonexpansiveness of \(p_\rho\) give \(L_p=\operatorname{Lip}(Df)/\rho\) as a Lipschitz bound for \(Dp_\rho\).

For adjacent cells, the difference quotient of the integrand is bounded by \(L_p\bigl((1-t)|\delta_{h,j}y_\ast|+t|\delta_{h,j}y_k|\bigr)\). Integrating the weights and taking the grid norm gives
\begin{equation}\label{eq:ss-average-jacobian-differences}
  \|\delta_h A_k\|_{h,\mathrm F}
   \le\tfrac12 L_p(\|\delta_h y_\ast\|_h+\|\delta_h y_k\|_h)
   \le C .
\end{equation}
The last bound uses the uniform \(H_h^1\) bounds on \(z_h^{(k)}\), \(z_h^\ast\) and the \(H_h^2\) bounds on \(x_h^{(k)}\), \(x_h^\ast\). In the exact error update \(e_{k+1}=(I-\alpha_k A_k)e_k+\alpha_k(2A_k-I)d_k\), let \(+\) denote the shift by \(h_je_j\). Taking its forward difference gives
\[
 \begin{aligned}
 \delta_{h,j}e_{k+1}
   &=(I-\alpha_k A_k^+)\delta_{h,j}e_k
      +\alpha_k(2A_k^+-I)\delta_{h,j}d_k\\
   &\quad+\alpha_k(\delta_{h,j}A_k)(2d_k-e_k).
 \end{aligned}
\]
The shifted matrices have the same spectral bounds. Summing over directions yields
\[
 \|\delta_h e_{k+1}\|_h
    \le\gamma_p\|\delta_h e_k\|_h
       +C\bigl(\|\delta_h d_k\|_h+\|d_k\|_{\infty,h}+\|e_k\|_{\infty,h}\bigr).
\]
The forcing decays with \(q_{\mathrm{DR}}\) by \eqref{eq:ss-recovery-smoothing}, the pointwise error bound, and \eqref{eq:ss-average-jacobian-differences}. Geometric summation from the bounded \(H_h^1\) error at \(k=1\) gives \(\|e_k\|_{H_h^1}\le Cq_{\mathrm{DR}}^{k-1}\).

For the nonlinear term in the residual, write
\[
 D_k=\int_0^1 Df((1-t)x_h^\ast+t x_h^{(k)})\,dt,\qquad
 f(x_h^{(k)})-f(x_h^\ast)=D_kd_k .
\]
Pointwise, \(\|D_k\|\le L_\star\). The same cellwise difference estimate as in \eqref{eq:ss-average-jacobian-differences} gives
\[
 \|\delta_h D_k\|_{h,\mathrm F}
   \le\frac{\operatorname{Lip}(Df)}{2}
      \bigl(\|\delta_h x_h^\ast\|_h+\|\delta_h x_h^{(k)}\|_h\bigr)\le C .
\]
The bound uses the uniform \(H_h^2\) norm of the recovered slopes. For each direction, the shifted product formula is
\[
  \delta_{h,j}(D_kd_k)
   =D_k^+\,\delta_{h,j}d_k+(\delta_{h,j}D_k)d_k .
\]
The first term is bounded by \(L_\star\|\delta_{h,j}d_k\|_h\) and the second by \(\|\delta_{h,j}D_k\|_{h,\mathrm F}\|d_k\|_{\infty,h}\). Combining with the grid norm bound on \(D_kd_k\) yields
\[
 \|f(x_h^{(k)})-f(x_h^\ast)\|_{H_h^1}
 \le C\bigl(\|d_k\|_{H_h^1}+\|d_k\|_{\infty,h}\bigr).
\]
Furthermore, \(\lambda_h^{(k)}-\lambda_h^\ast=\rho(d_k-e_k)\), so
\[
 \|f(x_h^{(k)})-\lambda_h^{(k)}\|_{H_h^1}
 \le C\bigl(\|d_k\|_{H_h^1}+\|d_k\|_{\infty,h}+\|e_k\|_{H_h^1}\bigr)
 \le Cq_{\mathrm{DR}}^{k-1}.
\]
Finally \(B_{T,h}^\ast=B_h^\ast\) in SS, with symbol norm \(\sqrt\nu\). Thus \(\|B_h^\ast w_h\|_h\le\|\delta_hw_h\|_h\), which proves the unnormalized bound from \eqref{eq:stationary-residual-identity}. In both models the constants depend on the box, model parameters, \(\rho\), \(\tau_0\), the relaxation bounds, and the uniform history and first-sweep bounds in the theorem.
\end{proof}

For a positive tolerance \(\mathrm{tol}\), it suffices to take
\[
 k=1+\max\biggl\{0,\biggl\lceil
    \frac{\log(C_R/\mathrm{tol})}{\log(1/q_{\mathrm{DR}})}
  \biggr\rceil\biggr\}.
\]
This bounds the number of inner sweeps for the residual test. The cost per sweep and the finite-precision acceptance checks are separate. Standard ADMM takes \(\alpha_k=1\). The estimate also applies to the accepted ADMM-Search steps in the prescribed interval. For PFC the projection removes the constant mode of the optimality residual. The time-equation defect multiplies each remaining residual mode by \(\nu_{\boldsymbol k}\). Applying the symbol bound of \(K_h\) to \eqref{eq:mesh-uniform-residual-bound} alone gives
\[
    \|K_hr_h^{(k)}\|_h
      \le 4\Big(\sum_{j=1}^d h_j^{-2}\Big) C_R q_{\mathrm{DR}}^{k-1}.
\]
A mesh-uniform estimate for the time-equation defect would therefore require additional control of the residual's spatial differences.

\paragraph{An a posteriori bound for the time-step solve}
Fix a discrete time step for SS or standard PFC and its history load. Let \(u_h\in U_h\) be an approximate solution, \(u_h^\ast\) the exact minimizer, and \(e_h=u_h-u_h^\ast\in T_h\). In PFC both fields have the prescribed mean. Write \(r_h(u_h)\) for the unnormalized projected gradient in \eqref{eq:normalized-optimality-residual} evaluated at \(u_h\). Since \(P_he_h=e_h\) and \(r_h(u_h^\ast)=0\), monotonicity of \(f\) and \eqref{eq:uniform-spectral-lower-bound} give
\[
  \begin{aligned}
    \langle r_h(u_h),e_h\rangle_h
     &= \langle Q_he_h,e_h\rangle_h
        +\langle f(B_hu_h)-f(B_hu_h^\ast),B_he_h\rangle_h\\
     &\ge c_0\|e_h\|_{H_h^2}^2 .
  \end{aligned}
\]
The spectral definition gives \(\langle r_h(u_h),e_h\rangle_h\le
\|r_h(u_h)\|_{H_h^{-2}}\|e_h\|_{H_h^2}\) by Cauchy--Schwarz with the factors \((I_h+K_h)^{-1}\) and \(I_h+K_h\), respectively. Dividing when \(e_h\ne0\) and using \((1+\nu_{\boldsymbol k})^{-1}\le1\) gives
\[
    \|u_h-u_h^\ast\|_{H_h^2}
      \le c_0^{-1}\|r_h(u_h)\|_{H_h^{-2}}
      \le c_0^{-1}\|r_h(u_h)\|_h .
\]
The bound controls the algebraic error for this fixed step. Discretization error relative to the continuous PDE and error propagation along a finite-tolerance trajectory require separate estimates.

\subsection{Exact outer trajectories}
\label{sec:mesh-uniform-outer-trajectories}

A restart from \(u_h^{(0)}=u_h^n\) uses either \(\lambda_h^{(0)}=0\) or \(\lambda_h^{(0)}=f(B_hu_h^n)\). The first sweep gives
\[
 z_h^{(1)}=v_h^{(1)}-\lambda_h^{(0)}/\rho,\qquad
 v_h^{(1)}=p_\rho\bigl(B_hu_h^n+\lambda_h^{(0)}/\rho\bigr).
\]
With the zero multiplier, nonexpansiveness and \(p_\rho(0)=0\) bound the value and first differences of \(z_h^{(1)}\) by those of \(B_hu_h^n\). With the nonlinear-gradient restart, \(v_h^{(1)}=B_hu_h^n\) by the proximal equation. If \(B_hu_h^n\) is uniformly bounded in \(H_h^2\), then the Lipschitz bound of \(f\) on its range also bounds \(\|z_h^{(1)}\|_{\infty,h}\) and \(\|\delta_h z_h^{(1)}\|_h\) for this restart. The \(H_h^2\) bound for \(B_hu_h^n\) follows from the single-step history assumptions for each model.

\begin{corollary}\label{cor:mesh-uniform-outer-trajectories}
Consider unforced BDF2 trajectories with BE startup and \(\sigma\ge1/16\), taking the exact minimizer at each outer time step. Assume the initial physical energies are uniformly bounded above and the initial SS means are uniformly bounded. Use the above restarts except at SS startup, where the guess is the sampled initial mean with zero multiplier. The ADMM inner iteration at each step retains the exact subproblem solves and the penalty and relaxation assumptions of Theorem~\ref{thm:mesh-uniform-admm}. Its bounds hold uniformly in the mesh and the time level, including the BE and first BDF2 solves. For PFC they are uniform in \(0<\tau\le\tau_0\), while for SS they hold at each fixed \(\tau>0\).
\end{corollary}

\begin{proof}
For SS, the starting energy and mean control \(u_h^0\) in \(H_h^2\). At the BE startup and first BDF2 step, we verify the bounds on the quadratic minimizer and the first-sweep state directly. Later steps satisfy the full \(H_h^3\) history assumption.

For \(a_h=u_h^1-u_h^0\), the corrections \(\frac{1}{4\tau}\|a_h\|_{L_h^\dagger,h}^2\) and \(\frac12\|a_h\|_{C_h,h}^2\) in \eqref{eq:modified-bdf2-energy} for \(\widetilde E_h^1\) are covered by the terms in the BE estimate \eqref{eq:be-energy-dissipation}. Thus \(\widetilde E_h^1\le E_h(u_h^0)\). The BDF2 energy inequality from \eqref{eq:bdf2-energy-estimate} gives
\(E_h(u_h^n)\le\widetilde E_h^n\le\widetilde E_h^1\) for \(n\ge1\).

For PFC, completing the square in the local potential gives
\[
 E_h(u_h)\ge
  \tfrac12\|(I_h-K_h)u_h\|_h^2
   +\tfrac18\|u_h\|_{4,h}^4-\tfrac12\epsilon^2|\Omega| .
\]
The quartic term controls \(\|u_h\|_h\), and the Fourier inequality \(\nu^2\le2(1-\nu)^2+2\) then controls \(\|K_hu_h\|_h\). This gives the required \(H_h^2\) histories and first-sweep bounds uniformly in the time level, the mesh, and \(0<\tau\le\tau_0\).

For SS the constant energy offset in \eqref{eq:mbe-ss-discrete-energy} transfers the bound to the physical energy, which controls \(\|K_hu_h^n\|_h\). The conserved mean and the positive lower bound on the nonzero symbols give a uniform \(H_h^2\) bound. Consequently \(\|f(B_hu_h^n)\|_h=\|B_hu_h^n\|_{6,h}^3\le C\) by \eqref{eq:mesh-discrete-embeddings} and the slope symbol estimate. We now fix \(\tau>0\).
At BE startup,
\[
 Q_{1,h}u_h^1=b_{1,h}^0-B_h^\ast f(B_hu_h^1),\qquad
 b_{1,h}^0=(\tau^{-1}I_h+K_h)u_h^0 .
\]
The load is bounded in \(L_h^2\) by the \(H_h^2\) energy bound, and the nonlinear term is bounded in \(H_h^{-1}\). The inverse of \(Q_{1,h}\) maps \(H_h^{-1}\) to \(H_h^3\), so \(u_h^1\) and the quadratic minimizer \(Q_{1,h}^{-1}b_{1,h}^0\) are bounded in \(H_h^3\). The constant SS guess has zero lifted slope and, with zero multiplier, gives \(z_h^{(1)}=0\).

For BDF2 write the exact equation as
\[
 Q_{2,h}u_h^{n+1}=\sigma\tau K_h^2u_h^n+c_h^n .
\]
The term \(c_h^n\) contains the first two terms of the load in \eqref{eq:bdf2-step-coefficients}. It also contains \(-B_h^\ast f(B_hu_h^{n+1})\). Write \(M_2\) for the uniform \(H_h^2\) bound from energy and mean. At fixed \(\tau\), the bounds on \(K_h:H_h^2\to L_h^2\) and \(B_h^\ast:L_h^2\to H_h^{-1}\) give
\[
 \|c_h^n\|_{H_h^{-1}}\le
   \bigl(5/(2\tau)+3\bigr)M_2+C_\Omega^3M_2^3 .
\]
The factor \(\theta_\tau=\sigma\tau/(\varepsilon_m^2+\sigma\tau)<1\) bounds the stabilization symbol divided by \(\widehat Q_{2,h}(\nu)\). As the quadratic inverse maps \(H_h^{-1}\) to \(H_h^3\) with norm at most \(1/c_0\), we have
\[
 \|u_h^{n+1}\|_{H_h^3}\le\theta_\tau\|u_h^n\|_{H_h^3}+C_\tau,\qquad n\ge1,
\]
where \(C_\tau\) is \(c_0^{-1}\) times the preceding bound on \(\|c_h^n\|_{H_h^{-1}}\). Iterating from the BE value gives, for \(n\ge1\),
\[
 \|u_h^n\|_{H_h^3}
 \le \theta_\tau^{n-1}\|u_h^1\|_{H_h^3}
     +C_\tau\frac{1-\theta_\tau^{n-1}}{1-\theta_\tau}.
\]
The right side is uniformly bounded in \(h\) and \(n\) at fixed \(\tau\).

At the first BDF2 solve, the history \(u_h^0\) only enters through \(u_h^0\) and \(K_hu_h^0\), both bounded in \(L_h^2\). The stabilizing load \(\sigma\tau K_h^2u_h^1\) is bounded in \(H_h^{-1}\) since \(u_h^1\) is bounded in \(H_h^3\). Thus the quadratic minimizer \(Q_{2,h}^{-1}b_{2,h}^1\) is already bounded in \(H_h^3\), and the restart from \(u_h^1\) gives the first-sweep bounds. These bounds furnish the inputs used in the state and residual proofs. At every subsequent solve both history fields are bounded in \(H_h^3\) by the recurrence.
\end{proof}

The corollary concerns exact outer solutions. Finite-tolerance evolution includes the signed residual work from the energy estimate, and reusing a previous multiplier requires its own verification of the first-sweep bounds.

\section{Numerical verification}
\label{sec:numerical-results}

\subsection{Smooth-data refinement test}

We test the first BDF2 solve on periodic \(N\times N\) cell-centered grids of \((0,2\pi)^2\), with \(h=2\pi/N\). We sample the smooth profiles
\[
    \begin{aligned}
        u_{0,\mathrm{SS}}(x,y)&=0.1[\sin(3x)\sin(2y)+\sin(5x)\sin(5y)],\\
        u_{0,\mathrm{PFC}}(x,y)&=-0.35+0.1\cos(x)\cos(y)
    \end{aligned}
\]
at cell centers and prescribe the sampled mean for PFC. We set \(\varepsilon_m=0.03\) for SS and \(\epsilon=0.35\) for PFC, with unit mobility, \(\sigma=1/8\), and \(\rho=1\). We vary \(N\) at each fixed step size listed in Table~\ref{tab:mesh-uniform-sweeps}. All cases are unforced.

For each case, BE with standard ADMM and tolerance \(10^{-10}\) supplies the common history \(u_h^1\). Both BDF2 solves use \((u_h^0,u_h^1)\), start from \(u_h^1\) and zero multiplier, and take an unrelaxed first sweep (\(\alpha_0=1\)). ADMM-Search uses \([0.1,1.9]\) thereafter. The finite-tolerance startup differs from the exact outer trajectory in Corollary~\ref{cor:mesh-uniform-outer-trajectories}.

The initial normalized residuals range from \(0.295\) to \(1\) over the 30 paired problems. For PFC the recorded initial residual is below one, so the normalization denominator is one. For SS the iteration CSVs contain only normalized values, which equal one to roundoff at initialization. For \(\tau=0.01\), \suppsection{sec:e8-reproducibility}{Reproducing the mesh and time-step study} reconstructs the unnormalized SS residual histories on three grids. None of these histories reaches the absolute tolerance \(10^{-8}\) before the normalized stopping rule terminates the solve.

All 60 solves met the internal threshold within 160 sweeps (Table~\ref{tab:mesh-uniform-sweeps}). At fixed \(\tau\), the PFC counts are constant across grids. From \(N=32\) to \(512\), SS adds at most three ADMM sweeps and five ADMM-Search sweeps. ADMM-Search uses fewer sweeps in every pair. The PFC residual histories nearly coincide across grids. The SS histories separate modestly (Fig.~\ref{fig:mesh-uniform-residual}).

\vbox{\relax
\begin{table}[H]
    \caption{Sweeps at the first BDF2 level, excluding BE startup. Each entry gives ADMM / ADMM-Search at the first verified \(R_h^{(k)}\le10^{-8}\) in the maintained Fourier representation.}
    \label{tab:mesh-uniform-sweeps}
    \centering
    \begin{tabular}{@{}lc*{5}{c}@{}}
    \toprule
    & & \multicolumn{5}{c}{Grid size \(N\)}\\
    \cmidrule(l){3-7}
    Model & \(\tau\) & 32 & 64 & 128 & 256 & 512\\
    \midrule
    SS & 0.04   & 48 / 25 & 50 / 26 & 50 / 28 & 51 / 28 & 51 / 30\\
    SS & 0.01   & 36 / 19 & 37 / 20 & 37 / 21 & 37 / 22 & 37 / 23\\
    SS & 0.0025 & 30 / 17 & 33 / 18 & 33 / 18 & 33 / 18 & 33 / 19\\
    \addlinespace
    PFC & 0.5     & 19 / 8 & 19 / 8 & 19 / 8 & 19 / 8 & 19 / 8\\
    PFC & 0.125   & 17 / 7 & 17 / 7 & 17 / 7 & 17 / 7 & 17 / 7\\
    PFC & 0.03125 & 17 / 8 & 17 / 8 & 17 / 8 & 17 / 8 & 17 / 8\\
    \bottomrule
\end{tabular}

\end{table}

\begin{figure}[H]
    \centering
    \includegraphics[width=.97\linewidth]{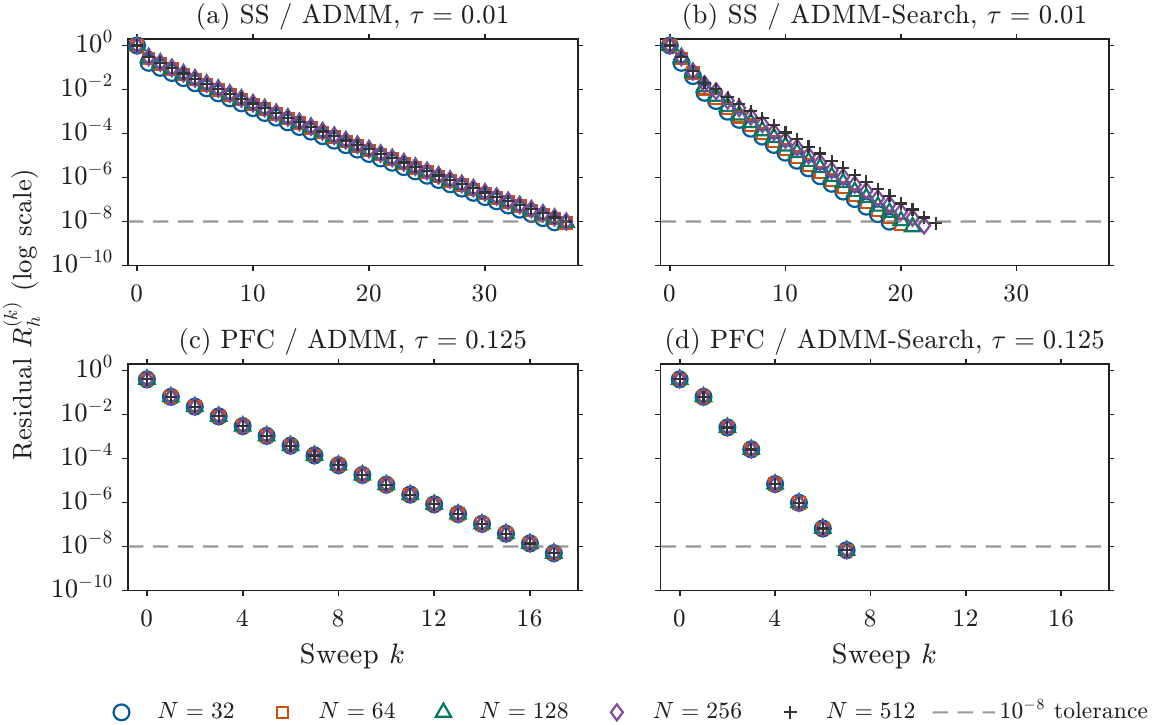}
    \caption{Normalized stopping residuals at the first BDF2 level for SS with \(\tau=0.01\) and PFC with \(\tau=0.125\). Panels (a) and (c) show ADMM, and panels (b) and (d) show ADMM-Search. The quadratic contribution uses the maintained Fourier coefficients and the nonlinear contribution is evaluated on the grid. Dashed lines mark the internal \(10^{-8}\) threshold.}
    \label{fig:mesh-uniform-residual}
\end{figure}
}

\vbox{\relax
\subsection{Stopping and output representations}

The \(10^{-8}\) threshold applies to the normalized stopping residual described in Fig.~\ref{fig:mesh-uniform-residual}. At acceptance the solver checks the time-equation defect and the physical energy for finiteness. The separate energy consistency check includes signed residual work, which can permit departures from energy monotonicity at finite tolerance.

On the finest PFC grid, reevaluating the residual on the returned grid arrays after rounding gives \(2.1\times10^{-7}\) to \(4.0\times10^{-7}\) in the original runs. For ADMM-Search at \(N=512\) and \(\tau=0.03125\), an independent higher-precision evaluation gives \(1.21\times10^{-9}\) for the maintained spectrum captured in a diagnostic replay and \(3.93\times10^{-7}\) for the returned grid array. The normalization denominator is one, and the residual reevaluated on the returned array exceeds the internal threshold. At high frequencies the quadratic PFC symbol grows as the square of the discrete Laplacian symbol, amplifying differences between representations. At the next step, the two-dimensional driver assembles the load using new transforms of the accepted grid histories.

For this single solve, \suppsection{sec:e8-reproducibility}{Reproducing the mesh and time-step study} reports unnormalized \(H_h^{-2}\) residuals of \(1.699\times10^{-10}\) for the returned array and \(1.696\times10^{-10}\) for the maintained spectrum. These norms are relevant to the \(H_h^2\) algebraic error estimate in Section~\ref{sec:mesh-uniform-stopping-residual}, which applies to fields with the prescribed mean. See \suppsection{sec:numerical-solver-efficiency}{Solver efficiency and sensitivity} for the 100-step comparisons and parameter scans.
\par}

\section{Conclusion}
\label{sec:conclusion}

For SS and standard PFC convex-splitting steps, Q-linear contraction of the DR state yields a mesh-uniform geometric bound on the unnormalized projected optimality residual in the unpreconditioned grid \(L^2\) norm. The SS transfer uses a first-difference bound before applying the slope adjoint. The single-step estimate is uniform for \(0<\tau\le\tau_0\) under the stated history and initialization assumptions and exact subproblem solves. The prescribed restarts along exact trajectories yield uniform bounds in the mesh and time level, with step-size uniformity for PFC and fixed-step bounds for SS. The smooth-data tests show sweep counts consistent with the mesh-uniform bound. Extending the trajectory estimates to finite-tolerance advancement with multiplier reuse requires control of residual work and restart regularity.

\section*{Acknowledgements}
OpenAI's ChatGPT and Codex assisted in drafting and revising the manuscript, developing and checking mathematical arguments, reviewing literature and citations, and writing and testing numerical and post-processing code. The authors assume responsibility for all content.

\bibliographystyle{siamplain}
\bibliography{references}

\end{document}